\documentclass[11pt]{article}
\usepackage[letterpaper,hmargin=1in,vmargin=1.25in]{geometry}

\usepackage{latexsym, amsfonts, amsmath, amssymb, amsthm}

\def\~{\tilde }

\def\ZZ{{\mathbb Z}}
\def\EE{{\mathbb E}}

\def\P{{\sf P}}
\def\Q{{\sf Q}}

\def\1{{\sf 1}}
\newcommand{\charf}[1]{\mbox{\raise.48ex\hbox{$\chi$}$_{#1}$}}
\def\bone#1{{\sf 1}_{#1}}

\usepackage{graphicx}
\def\overloop#1{{\overset{\rotatebox{180}{$\circlearrowright$}}{#1}}}

\usepackage{semtrans}
\newcommand{\overphi}{\raise1.6ex\hbox{\rotatebox{180}{\textup{\U{$\Phi$}}}}}

\theoremstyle{plain}
\newtheorem{theorem}{Theorem}[section]
\newtheorem{lemma}[theorem]{Lemma}

\newtheorem{definition}[theorem]{Definition}

\theoremstyle{definition}
\newtheorem{example}[theorem]{Example}
\newtheorem{remark}[theorem]{Remark}

\begin{document}

\title{Intersection conductance and canonical alternating paths, \\
    methods for general finite Markov chains}
\author {Ravi Montenegro \thanks{Department of Mathematical Sciences, University of
   Massachusetts Lowell, Lowell, MA 01854, ravi\_montenegro@uml.edu.}
}

\maketitle
 
\begin{abstract}\noindent
We extend the conductance and canonical paths methods to the setting of general finite Markov chains, including non-reversible non-lazy walks.
The new path method is used to show that a known bound for mixing time of a lazy walk on a Cayley graph with symmetric generating set also applies to the non-lazy non-symmetric case, often even when there is no holding probability.

\vspace{1ex}
\noindent {\bf Keywords} : Mixing time, conductance, canonical paths, evolving sets, Cayley graph.
\end{abstract}

\section{Introduction}

Beginning with the work of Jerrum and Sinclair, geometric concepts such as conductance \cite{JS88.1, LS88.1} and canonical paths \cite{Sin92.1,DS91.1} have played an important role in studying the mixing time of finite ergodic Markov chains.
These methods originally applied only to reversible lazy walks, and while little is lost in dropping reversibility and some amount of laziness \cite{DS91.1,Fill91.1,MT06.1}, extensions which allow for dropping both conditions tend to be weak or difficult to use \cite{Fill91.1,Mih89.1,MT06.1}.
Similar difficulty has been encountered with the method of blocking conductance \cite{KLM06.1}, a geometric approach to sharpening conductance bounds by including a notion of vertex congestion but which has only been successfully applied to a few problems.
In this paper we develop an extension of conductance to the general (non-lazy non-reversible) setting, along with a further extension which sharpens bounds by using a simple notion of vertex congestion.
This makes it easier to use and allows for the proof of a new canonical path theorem which applies to general finite ergodic Markov chains.

Recall that if a finite Markov kernel $\P$ with sample space $V$ is irreducible and aperiodic ($\exists N\in{\mathbb N},\,\forall x,y\in V:\,\P^N(x,y)>0$) then it has a unique stationary distribution $\pi$ satisfying $\pi\P=\pi$ and is ergodic ($\forall x\in V:\,\P^n(x,\cdot)\xrightarrow{n\to\infty} \pi$).
The walk is lazy if the holding probability $\alpha=\min_{v\in V} \P(v,v)$ is at least $1/2$, and is reversible if the time-reversal (adjoint) $\P^*(x,y)=\frac{\pi(y)\P(y,x)}{\pi(x)}$ satisfies $\P^*=\P$.
The $L^2$ mixing time $\tau(\epsilon)=\max_{x\in V}\min\{n\,:\,\|k_n^x-1\|_{2,\pi}\leq\epsilon\}$ is the worst-case time for standard deviation of the density $k_n^x(y)=\frac{\P^n(x,y)}{\pi(y)}$ to drop to $\epsilon$.

Jerrum and Sinclair \cite{JS88.1} and Lawler and Sokal \cite{LS88.1} showed that mixing time of a lazy reversible walk can be bounded in terms of the conductance, also known as the Cheeger constant,
\[
\Phi=\min_{\emptyset \subsetneq A\subsetneq V} \Phi(A)
\qquad\textrm{where}\qquad
\Phi(A)=\frac{\Q(A,A^c)}{\pi(A)\pi(A^c)}
\]
and the ergodic flow from $A\subset V$ to $B\subset V$ is $\Q(A,B)=\sum_{x\in A,\,y\in B} \pi(x)\P(x,y)$.

\begin{theorem} \label{thm:conductance}
The mixing time of a lazy reversible finite ergodic Markov chain is
\[
\tau(\epsilon) 
  \leq \frac{8}{\Phi^2}\log \frac{1}{\epsilon\sqrt{\pi_0}}
\]
where $\pi_0=\min_{v\in V}\pi(v)$.
\end{theorem}

Since $\pi$ is stationary then $\pi(v)=\Q(A,v)+\Q(A^c,v)$ and $\Q(A,A^c)=\Q(A^c,A)$, and so in the case of a lazy walk 
\[
\frac 12\,\sum_{v\in V} \min\{\Q(A,v),\,\Q(A^c,v)\} = \frac{\Q(A,A^c)+\Q(A^c,A)}{2} = \Q(A,A^c)
\]
i.e. the flow $\Q(A,\cdot)$ from $A$ and $\Q(A^c,\cdot)$ from $A^c$ intersect at exactly the flow $\Q(A,A^c)+\Q(A^c,A)$ crossing between $A$ and $A^c$.
The intersection of the flows from $A$ and $A^c$ can then be used to define a modified form of conductance which agrees with normal conductance in the lazy case.

\begin{definition}
The {\em intersection conductance $\overphi(A)$} of $A\subset V$ (and $A^c$) is given by
\[
\overphi(A) = \overphi(A,A^c) = \frac{\sum_{v\in V} \min\{\Q(A,v),\,\Q(A^c,v)\}}{2\pi(A)\pi(A^c)}
\]
The {\em intersection conductance} is $\displaystyle \overphi = \min_{\emptyset \subsetneq A \subsetneq V} \overphi(A)$.
\end{definition}

Intuitively, when flow from $A$ and $A^c$ overlap significantly then a walk starting in $A$ will quickly mix with the `unoccupied' space starting in $A^c$, and so good mixing may be expected.
As our first result, in Section \ref{sec:threshold-conductances} we make this rigorous by showing that mixing time for general (non-lazy non-reversible) finite Markov chains can be bounded using the intersection conductance.

\begin{theorem} \label{thm:intersection-conductance}
The mixing time of a finite ergodic Markov chain is
\[
\tau(\epsilon) 
  \leq \frac{12}{\overphi^2}\log \frac{1}{\epsilon\sqrt{\pi_0}}
\]
\end{theorem}

As a further extension, suppose that a threshold $t$ is fixed and ergodic flow is counted up to at most a $t$ fraction of vertex capacity, 
i.e. work with {\em threshold limited ergodic flow} $\Q_t(A,v) = \min\{\Q(A,v),\,t\pi(v)\}$ instead of $\Q(A,v)$.
When ergodic flow is well distributed among vertices then it may be possible to make $t$ quite small without decreasing the ergodic flow significantly, and so the optimal choice of $t$ will measure some form of vertex congestion.
More formally, and in the greater generality of conductance profile:

\begin{definition} \label{def:intersection-threshold-conductance}
Given $t > 0$ the {\em intersection threshold conductance} of $A\subset V$ is
\[
\overphi_t(A) = \frac{\sum_{v\in V} \min\{\Q_t(A,v),\,\Q_t(A^c,v)\} }{2\pi(A)\pi(A^c)}
\]
The {\em intersection threshold conductance profile $\overphi_t(r)$} is 
\[
\overphi_t(r)=
\begin{cases}
\min_{0<\pi(A)\leq r}\overphi_t(A) & \textrm{if } r\leq 1/2 \\
\overphi_t(1/2) & \textrm{if } r>1/2
\end{cases}
\]
The {\em intersection threshold conductance} is $\displaystyle \overphi_t = \min_{\emptyset \subsetneq A \subsetneq V} \overphi_t(A) = \overphi_t(1/2)$.
\end{definition}

For a lazy walk this is equal to normal conductance when $t\geq 1/2$, i.e. $\overphi_{1/2}(A) = \Phi(A)$, while more generally $\overphi_{1/2}(A)=\overphi(A)\leq 1$.
Our earlier result can be sharpened to show that the mixing time bound is $t$ times smaller when a threshold is used.
In particular:

\begin{theorem} \label{thm:intersection-threshold-conductance}
Given threshold $t>0$ the mixing time of a finite ergodic Markov chain is
\[
\tau(\epsilon) \leq \frac{12\max\{t,\overphi_t\}}{\overphi_t^2}\log \frac{1}{\epsilon\sqrt{\pi_0}}
\]
and more generally
\[
\tau(\epsilon) \leq
\begin{cases}
\bigskip
\displaystyle \int_{4\pi_0}^{4/\epsilon^2} \frac{12\,\max\{t,\overphi_t(r)\}}{r\overphi_t(r)^2}\,dr & \textrm{in general} \\
\displaystyle \int_{\pi_0}^{1/\epsilon^2} \frac{6\,\max\{t,\overphi_t(r)\}}{r\overphi_t(r)^2}\,dr\hspace{0.25in} & \textrm{if $r\,\frac{\overphi_t^2\left(\frac{1}{1+r^2}\right)}{\max\left\{t,\overphi_t\left(\frac{1}{1+r^2}\right)\right\}}$ is convex in $r$}
\end{cases}
\]
\end{theorem}


In their seminal paper Jerrum and Sinclair introduced the method of canonical paths as a means of lower bounding conductance \cite{JS88.1}.

\begin{definition} \label{def:canonical-paths}
A {\em canonical path} $\gamma_{xy}$ is a path from $x$ to $y$ using only valid transitions of $\P$:
\[
x=x_0\rightarrow x_1 \rightarrow x_2 \rightarrow x_3 \rightarrow \cdots \rightarrow x_{n-1} \rightarrow x_n=y
\]

Let $\Gamma=\{\gamma_{xy}:\,x,y\in V\}$ include a canonical path for each pair of distinct vertices $x,y\in V$.

The {\em edge-congestion} is given by
\[
\rho_e = \rho_e(\Gamma) 
 = \max_{u,v\in V:\,\P(u,v)>0} \frac{1}{\pi(u)\P(u,v)}\,\sum_{(u,v)\in\gamma_{xy}\in\Gamma} \pi(x)\pi(y)
\]
\end{definition}


\begin{theorem}[Jerrum and Sinclair] \label{thm:JS}
The mixing time of a lazy reversible finite ergodic Markov chain is
\[
\tau(\epsilon) \leq 8\rho_e^2\log\frac{1}{\epsilon\sqrt{\pi_0}}
\]
\end{theorem}

Sinclair later showed a Poincar\'e type bound with order $\rho_e^2$ replaced by $\rho_e\ell$ where $\ell=\max |\gamma_{xy}|$ is the length of the longest path, usually an improvement on the $\rho_e^2$ result \cite{Sin92.1}.
Diaconis and Strook extend this to non-lazy reversible walks if $\Gamma$ also includes odd length paths $\gamma_{xx}$ from each vertex to itself \cite{DS91.1}.
Various authors have observed that the Poincar\'e approach applies to lazy non-reversible walks as well, and to non-lazy walks at the cost of a factor of $\alpha^{-1}$ 
\cite{Fill91.1,Mih89.1,MT06.1}.

As our second main result, in Section \ref{sec:alternating-paths} we find that a generalization of the canonical path method can be used to lower bound intersection threshold conductance, giving a canonical path method for general finite Markov chains.
In particular we consider {\em canonical alternating paths}, even length paths which alternate between forward and reversed edges of $\P$, i.e.
\[
x=x_0\rightarrow x_1 \leftarrow x_2 \rightarrow x_3 \leftarrow \cdots \rightarrow x_{2n-1} \leftarrow x_{2n}=y
\]
This is equivalent to a path from $x$ to $y$ alternating between edges of $\P$ and $\P^*$.
The canonical paths used in each of the methods discussed in the previous paragraph can be used to construct canonical alternating paths -- for instance by adding self-loops at $x_1$, $x_2$, $\ldots$, $x_n$ in the lazy case -- so this new type of path generalizes each of those settings.
For an appropriately defined notion of vertex congestion $\rho_v$ we show that
\[
\overphi_{\rho_v/\rho_e} \geq 1/2\rho_e\,.
\]
A mixing time result then follows from Theorem \ref{thm:intersection-threshold-conductance}.
With some extra work a bound is also possible if standard canonical paths are used:
\begin{theorem}
A finite ergodic Markov chain with canonical alternating paths $\Gamma$ has mixing time
\[
\tau(\epsilon) \leq 48 \rho_e \max\{\rho_v,\,\frac 12\} \log\frac{1}{\epsilon\sqrt{\pi_0}} 
\]
If ordinary canonical paths are used and $\alpha=\min_{v\in V}\P(v,v)$ is the minimal holding probability then
\[
\tau(\epsilon) \leq 4\rho_v\,\max\left\{\frac{\rho_v}{\alpha},\,\rho_e\right\}\,\log\frac{1}{\epsilon\sqrt{\pi_0}}
\]
\end{theorem}
\noindent
The congestions satisfy $\rho_v\leq\rho_e$, so this generalizes Jerrum and Sinclair's result.
Our results also hold for multicommodity flows.

We finish in Section \ref{sec:examples} with a few examples of how our new tools can be used.
First, it is shown that known complexity results for the lazy simple random walk on an undirected graph apply just as well to the simple random walk on an Eulerian directed graph, i.e. a strongly connected graph with in-degree=out-degree at each vertex, both when there is a self-loop at each vertex and often even without self-loops.
A more interesting problem is to study random walks on Cayley graphs, i.e. walks on groups, for which we show that known bounds \cite{Bab91.1,DSC93.2} for a lazy walk or a walk with a symmetric set of generators can be extended to the non-lazy non-symmetric case.


\section{Evolving Sets} \label{sec:evosets}

The Evolving Set methodology \cite{MP05.1,MT06.1} will be required for the proof of Theorem \ref{thm:intersection-threshold-conductance}, from which most results of the paper follow.
We give here a very brief introduction to those elements required in our proof.
The reader interested only in the canonical path results may skip to Section \ref{sec:alternating-paths}.

One approach to relating a property of sets (conductance) to a property
of the original walk (mixing time) is to construct a {\em dual process}:
a walk $\P_D$ on some state space $\Omega$ and a {\em
link}, or transition matrix, $\Lambda$ from $\Omega$ to $V$ such that
$$
\P\Lambda = \Lambda\P_D\,.
$$
In particular, $\P^n\Lambda = \Lambda\P_D^n$ and so the evolution
of $\P^n$ and $\P_D^n$ will be closely related. 
A natural candidate to link a walk on sets to a walk on states is
the projection $\Lambda(S,y) = \frac{\pi(y)}{\pi(S)}\,\1_{S}(y)$.
Diaconis and Fill \cite{DF90.1} have shown that for certain
classes of Markov chains a walk based on the Evolving Set process discussed below is the
unique dual process with link $\Lambda$, so this is the most natural walk on
sets to consider. 
Our discussion of Evolving Sets will be based on work of Morris and
Peres \cite{MP05.1} in a slightly improved form of Montenegro and Tetali \cite{MT06.1}.

To understand the method requires some new terminology.

\begin{definition}\label{def:evosets}
Given set $A\subset V$ a step of the {\em evolving set process} is given by choosing $u\in[0,1]$ uniformly at random, and transitioning to the set 
\[
A_u = \{y\in V\,\mid\,\Q(A,y)\geq u\pi(y)\} = \{y\in V\,\mid\,\P^*(y,A)\geq u\}\,.
\]
The {\em root profile} $\psi:\,(0,\infty)\to[0,1]$ is given by
\[
\psi(r)=\min_{0<\pi(A)\leq r} \psi(A)
\quad\textrm{where}\quad
\psi(A)=1-\frac{\int_0^1 \sqrt{\pi(A_u)(1-\pi(A_u))}\,du}{\sqrt{\pi(A)(1-\pi(A))}}
\]
when $r\in(0,1/2]$, and $\psi(r)=\psi(1/2)=\min_{A\not\in\{\emptyset,V\}}\psi(A)$ when $r>1/2$.
\end{definition}

The claim $\psi(1/2)=\min_{A\not\in\{\emptyset,V\}}\psi(A)$ follows from the relation $\psi(A)=\psi(A^c)$, a consequence of $(A_u)^c=(A^c)_{1-u}$ for every $u$ with $\nexists v:\,\Q(A,v)=u\pi(v)$, i.e. $u$-a.e. since $V$ is finite.

Since $\Q(V,v)=\pi(v)$ then $A_u$ consists of those vertices receiving at least a $u$-fraction of their steady state probability from $A$.
If $u$ is chosen uniformly from $[0,1]$ then
\[
\EE \pi(A_u) = \int_0^1 \pi(A_u)\,du=\sum_{y\in V}\pi(y)\frac{\Q(A,y)}{\pi(y)}=\pi(A)
\]
and so by Jensen's Inequality $\EE \sqrt{\pi(A_u)(1-\pi(A_u))} \leq \sqrt{\pi(A)(1-\pi(A))}$, with equality if and only if $\pi(A_u)=\pi(A)$ $u$-a.e.
It follows that a large root profile $\psi(A)$ indicates that $\pi(A_u)$ differs significantly from $\pi(A)$,
and in particular the flow from $A$ is spread over a large space and so the walk expands quickly from $A$.
Morris and Peres \cite{MP05.1} make this intuition rigorous, although we use a sharper result of Montenegro and Tetali \cite{MT06.1}:

\begin{theorem} \label{thm:evomixing}
A finite ergodic Markov chain with root profile lower bounded by $\psi(r)$ has
\[
\tau(\epsilon) \leq 
\begin{cases}
\vspace{2ex}\displaystyle \int_{4\pi_0}^{4/\epsilon^2} \frac{dr}{r\,\psi(r)} 
 & in\ general \\
\displaystyle \int_{\pi_0}^{1/\epsilon^2} \frac{dr}{2r\psi(r)} 
 & if\ r\psi\left(\frac{1}{1+r^2}\right)\ is\ convex
\end{cases}
\]
\end{theorem}

The root profile is typically lower bounded by writing constraints of interest in terms of the evolving set process.
For instance, when conductance is being considered then a lazy walk has $A\subset A_u$ when $u\leq 1/2$ and $A_u\subset A$ when $u>1/2$ and so
\begin{equation*}
\Q(A,A^c) 
 = \sum_{v\in A^c} \frac{\Q(A,v)}{\pi(v)}\,\pi(v) 
 = \sum_{v\in A^c} \int_0^{1/2} \bone{\Q(A,v)/\pi(v)\geq u}\,\pi(v)\,du 
 = \int_0^{1/2} \pi(A_u\setminus A)\,du
\end{equation*}
The identities $\int_0^1 \pi(A_u)\,du = \pi(A)$ and $\Q(A,A^c)=\Q(A^c,A)$ can be used to write this as an area problem:
\begin{equation} \label{eqn:evolving-flow}
\Q(A,A^c) = \int_0^{1/2} (\pi(A_u)-\pi(A))\,du = \int_{1/2}^1 (\pi(A)-\pi(A_u))\,du = \Q(A^c,A)
\end{equation}
Breaking the definition of $\psi(A)$ into an integral of $u\in[0,1/2]$ and one of $u\in[1/2,1]$, applying Jensen's Inequality, and then making a few simplifications leads to the relation $\psi(A)\geq\Phi(A)^2/2$.
Theorem \ref{thm:evomixing} then implies a stronger version of Jerrum and Sinclair's result given in Theorem \ref{thm:conductance}:
\begin{theorem} \label{thm:conductance-improved}
A lazy finite ergodic Markov chain has $\psi(A) \geq \Phi^2(A)/2$ and mixing time
\[
\tau(\epsilon) \leq \frac{2}{\Phi^2}\,\log \frac{1}{\epsilon\sqrt{\pi_0}}
\]
\end{theorem}


\section{Mixing bounds with Threshold Conductances} \label{sec:threshold-conductances}

In this section we prove a lower bound on the root profile $\psi(A)$ in terms of the (intersection) threshold conductance of Definition \ref{def:intersection-threshold-conductance}. 
By Theorem \ref{thm:evomixing} this induces upper bounds on mixing time, including the main result of this paper, Theorem \ref{thm:intersection-threshold-conductance}, a bound on mixing time in terms of the intersection threshold conductance.
To further indicate the improvement provided by use of thresholds we also give a bound in terms of a quantity which more strongly resembles ordinary conductance:

\begin{definition}
If $A\subset V$ and $t\in[0,1]$ then the {\em threshold conductance $\hat\Phi_t(A)$} is
\[
\hat\Phi_t(A)=\frac{\min\{\Q_t(A,A^c),\,\Q_t(A^c,A)\}}{\pi(A)\pi(A^c)}
\]
The {\em threshold conductance profile $\hat\Phi_t(r)$} and the {\em threshold conductance $\hat\Phi_t$} are defined as were the related quantities in Definition \ref{def:intersection-threshold-conductance}.
\end{definition}

For a lazy walk $\Q_{1/2}(A,A^c)=\Q(A,A^c)$ and $\hat\Phi_{1/2}(A)=\Phi(A)$ agree with standard notions of ergodic flow and conductance.
Intuitively, if ergodic flow from $A$ to $A^c$ and from $A^c$ to $A$ are not unduly concentrated on a few vertices then $\hat\Phi_t(A)\approx\Phi(A)$ even for fairly small $t$.
The extra information provided by $t$ will be used here to show a substantially improved result.

And now to the main result:

\begin{theorem} \label{thm:threshold-conductances}
Given $A\subset V$ then 
\[
\psi(A) 
\geq
\begin{cases}
\bigskip\displaystyle\frac{1}{12t}\,\min\left\{\overphi_t(A)^2,\,t\overphi_t(A)\right\} & \\
\displaystyle\frac{\min\{\alpha,t\}}{4t^2}\,\hat\Phi_t(A)^2 & 
\end{cases}
\]
\end{theorem}


The second bound simply sharpens a special case of the first, as 
$\overphi_{\min\{\alpha,t\}}(A) \geq \frac{\min\{\alpha,t\}}{t}\,\hat\Phi_t(A) \leq 2\min\{\alpha,t\}$.
When $\hat\Phi_t(r)=\Theta(\Phi(r))$ then this leads to a mixing bound $\frac{t^2}{\min\{\alpha,t\}}$ times that of Theorem \ref{thm:conductance}, a significant improvement when the threshold is small.
In the extreme case when ergodic flow is spread uniformly over the complement then $\hat\Phi_t(r)=\Phi(r)$ when $t=\min \frac{\Q(A,A^c)}{\min\{\pi(A),\pi(A^c)\}}\leq\Phi$,
improving the mixing time bound from $\tau(\epsilon)\leq \frac{2}{\Phi^2}\log\frac{1}{\epsilon\sqrt{\pi_0}}$ to
$\tau(\epsilon) = O\left(\frac{1}{\Phi}\log\frac{1}{\epsilon\sqrt{\pi_0}}\right)$ if $\alpha\geq\Phi$.
With conductance profile this can even match the best case lower bound of $\tau(\epsilon)=\Omega\left(\frac{1}{\Phi}\right)$.

\begin{proof}
Start with the first bound of the theorem, that $\psi(A) \geq \frac{1}{12t}\,\min\left\{\check\phi_t(A)^2,\,t\check\phi_t(A)\right\}$.
The proof will involve relating evolving sets to the numerator of $\overphi_t(A)$, with \eqref{eqn:evolving-intersection} replacing the use of \eqref{eqn:evolving-flow} for conductance.
Breaking the definition of $\psi(A)$ into an integral over $u\in[0,t]$ and one over $u\in[t,1]$, applying Jensen's Inequality, and then making a few simplifications completes the proof.
The argument for the second case is similar.

Suppose $t\geq 1/2$.
Since $\min\left\{\Q(A,v),\,\Q(A^c,v)\right\} \leq \pi(v)/2$ then $\min\left\{\Q_t(A,v),\,\Q_t(A^c,v)\right\}=\min\left\{\Q_{1/2}(A,v),\,\Q_{1/2}(A^c,v)\right\}$ 
and consequently $\overphi_t(A)=\overphi_{1/2}(A)$. 
The result at $t=1/2$ is then stronger than that at $t>1/2$, so without loss assume that $t\leq 1/2$.

When $u\leq 1/2$ then
\begin{eqnarray*}
\pi(A_u\setminus A_{1-u})
 &=& \pi\left(\left\{v\in V:\,\Q(A,v)\geq u\pi(v),\,\Q(A^c,v)> u\pi(v)\right\}\right) \\
 &=& \pi\left(\left\{v\in V:\,\min\{\Q(A,v),\,\Q(A^c,v)\}\geq u\pi(v)\right\}\right) \\
 &&- \pi\left(\left\{v\in V:\,\Q(A,v)\geq u\pi(v),\,\Q(A^c,v)=u\pi(v)\right\}\right) 
\end{eqnarray*}
The set $\{u\in[0,1]:\,\exists v\in V,\,\Q(A^c,v)=u\pi(v)\}$ is finite and so it has Lebesgue measure zero.
Since $u\leq 1/2$ then $A_{1-u}\subset A_u$ and so
\begin{eqnarray}
\lefteqn{\sum_{v\in V} \min\left\{\Q_t(A,v),\,\Q_t(A^c,v)\right\} = \sum_{v\in V} \int_0^t \pi(v)\,\delta_{\min\left\{\Q_t(A,v),\,\Q_t(A^c,v)\right\}\geq u\pi(v)}\,du }  \nonumber \\
 &=& \int_0^t \left[\pi(A_u\setminus A_{1-u})
  + \pi\left(\left\{v\in V:\,\Q(A,v)\geq u\pi(v),\,\Q(A^c,v)=u\pi(v)\right\}\right)\right]\,du  \nonumber \\
 &=& \int_0^t \left(\pi(A_u)-\pi(A_{1-u})\right)\,du + 0  \nonumber \\
 &=& \int_0^t (\pi(A_u)-\pi(A))\,du + \int_{1-t}^1 (\pi(A)-\pi(A_u))\,du \label{eqn:evolving-intersection} \\
 &=& \int_0^t (\pi(A_u)-\pi(A))\,du + \int_0^t (\pi((A^c)_u)-\pi(A^c))\,du  \nonumber 
\end{eqnarray}
The final equality uses the relation $(A_u)^c = (A^c)_{1-u}$ for a.e. $u\in[0,1]$.
It follows that
\[
\max\left\{\int_0^t (\pi(A_u)-\pi(A))\,du,\, \int_0^t (\pi((A^c)_u)-\pi(A^c))\,du\right\} \geq \overphi_t(A)\pi(A)\pi(A^c)
\]
Since $\psi(A)=\psi(A^c)$ and $\overphi_t(A)=\overphi_t(A^c)$ then the roles of $A$ and $A^c$ can be swapped and this inequality will still hold.
So without loss assume $\int_0^t (\pi(A_u)-\pi(A))\,du \geq \int_0^t (\pi((A^c)_u)-\pi(A^c))\,du$. By Jensen's Inequality, the Martingale identity $\int_0^1 \pi(A_u)\,du=\pi(A)$, and concavity of $f(x)=\sqrt{x(1-x)}$ then
\begin{eqnarray*}
\lefteqn{\int_0^1 f(\pi(A_u))\,du
  \leq t\,f\left(\int_0^t \pi(A_u)\,\frac{du}{t}\right)+ (1-t)\,f\left(\int_t^1 \pi(A_u)\,\frac{du}{1-t}\right) } \\
  &=&    t\,f\left(\pi(A) + \frac{\int_0^t (\pi(A_u)-\pi(A))\,du}{t}\right) + (1-t)\,f\left(\pi(A) - \frac{\int_0^t (\pi(A_u)-\pi(A))\,du}{1-t}\right)
\end{eqnarray*}
By concavity of $f$ this is decreasing in $\int_0^t (\pi(A_u)-\pi(A))\,du\geq \overphi_t(A)\pi(A)\pi(A^c)$ and so it is maximized at the lower bound.
The root profile is then bounded by
\begin{eqnarray}
1-\psi(A) &\leq& t\sqrt{\left(1+\frac{\overphi_t(A)(1-\pi(A))}t\right)\left(1-\frac{\overphi_t(A)\pi(A)}t\right)} \label{eqn:Schur} \\
    && + (1-t)\sqrt{\left(1-\frac{\overphi_t(A)(1-\pi(A))}{1-t}\right)\left(1+\frac{\overphi_t(A)\pi(A)}{1-t}\right)} \nonumber
\end{eqnarray}
Treating $t$ and $\overphi_t(A)$ as constants, this is concave in $\pi(A)$ with negative derivative at $\pi(A)=0$ when $\overphi_t(A)<1-2t$.
The bound in \eqref{eqn:Schur} is then maximized at $\pi(A)=0$, with
\begin{eqnarray*}
1-\psi(A) &\leq& t\sqrt{1+\frac{\overphi_t(A)}t} + (1-t)\sqrt{1-\frac{\overphi_t(A)}{1-t}} \\
 &\leq& t\left(1+\frac{\overphi_t(A)}{2t}-\frac{1}{12}\min\left\{\frac{\overphi_t(A)^2}{t^2},\frac{\overphi_t(A)}t\right\}\right)
    + (1-t)\left(1-\frac{\overphi_t(A)}{2(1-t)}\right)
\end{eqnarray*}
The second inequality applies the relations $\sqrt{1+x}\leq 1+\frac{x}{2}-\delta_{x\geq 0}\frac{\min\{x^2,x\}}{12}$ and $\sqrt{1-x}\leq 1-\frac x2$.
This gives the theorem when $\overphi_t(A)<1-2t$.
When $\overphi_t(A)\geq 1-2t$ then \eqref{eqn:Schur} is maximized at $\pi(A)=\frac{1-(1-2t)/\overphi_t(A)}{2}$, leading to the bound $\psi(A) \geq \frac{\overphi_t(A)^2}{2} \geq \frac{1}{6t}\,\min\{\overphi_t(A)^2,\,t\overphi_t(A)\}$.

Now consider the second bound of the theorem, $\psi(A)\geq\frac{\min\{\alpha,t\}}{4t^2}\,\hat\Phi_t(A)^2$.

Suppose $t\geq 1/2$.
If $\alpha\geq 1/2$ then $\hat\Phi_t(A)=\hat\Phi_{1/2}(A)$ and so the bound follows from the case of $t=1/2$.
If $\alpha<1/2$ then $\hat\Phi_{1/2}(A)\geq\frac{1}{2t}\hat\Phi_t(A)$ and again the bound follows from the case of $t=1/2$.
The result at $t=1/2$ is then stronger than that at $t>1/2$, so without loss assume that $t\leq 1/2$.

Let $T=\min\{\alpha,t\}\leq 1/2$. Observe that $\Q_T(A,v) = \int_0^T \pi(A_u\cap \{v\})\,du$ and so 
\begin{eqnarray*}
\Q_T(A,V) &=& \sum_{v\in V} \int_0^T \pi(A_u\cap \{v\})\,du = \int_0^T \pi(A_u)\,du \\
\Q_T(A,V) &=& \Q_T(A,A) + \Q_T(A,A^c) = \pi(A)T + \Q_T(A,A^c)
\end{eqnarray*}
Then $\int_0^T \pi(A_u)\,du = \Q_T(A,V) = \pi(A)T + \Q_T(A,A^c)$.
Swapping the roles of $A$ and $A^c$ in this identity, and recalling that $u$-a.e. $(A_u)^c=(A^c)_{1-u}$, it follows that
\[
\int_{1-T}^1 \pi(A_u)\,du = \int_0^T (1-\pi((A^c)_u))\,du = \pi(A)T - \Q_T(A^c,A)
\]

By Jensen's Inequality and the identity $\int_0^1 \pi(A_u)\,du = \pi(A)$:
\begin{eqnarray*}
\lefteqn{\int_0^1 \sqrt{\pi(A_u)(1-\pi(A_u))}\,du 
  = \int_0^{T} + \int_{T}^{1-T} + \int_{1-T}^1 f(\pi(A_u))\,du } \\
  &\leq& T\,f\left(\int_0^{T}\pi(A_u)\,\frac{du}{T}\right) + (1-2T)\,f\left(\int_T^{1-T}\pi(A_u)\,\frac{du}{1-2T}\right)  + T\,f\left(\int_{1-T}^1\pi(A_u)\,\frac{du}{T}\right) \\
  &=& T\,f\left(\pi(A) + \frac{\Q_T(A,A^c)}{T}\right) + (1-2T)\,f\left(\pi(A) - \frac{\Q_T(A,A^c)-\Q_T(A^c,A)}{1-2T}\right) \\
  && + T\,f\left(\pi(A) - \frac{\Q_T(A^c,A)}{T}\right)
\end{eqnarray*}

This is decreasing in $\Q_T(A,A^c)$ when $\Q_T(A,A^c)\geq\Q_T(A^c,A)$, and decreasing in $\Q_T(A^c,A)$ when $\Q_T(A^c,A)\geq\Q_T(A,A^c)$, and so is maximized when $\Q_T(A,A^c)=\Q_T(A^c,A)=\hat\Phi_t(A)\pi(A)\pi(A^c)$.
It follows that
\begin{eqnarray*}
1-\psi(A) &\leq& \frac{T\,f\left(\pi(A) + \frac{\hat\Phi_t(A)\pi(A)\pi(A^c)}{T}\right) + (1-2T)\,f(\pi(A)) + T\,f\left(\pi(A) - \frac{\hat\Phi_t(A)\pi(A)\pi(A^c)}{T}\right)}{f(\pi(A))} \\
&=& 2T\,\sqrt{\left(\frac 12+\frac{\hat\Phi_T(A)}{2T}\pi(A^c)\right)\left(\frac 12-\frac{\hat\Phi_T(A)}{2T}\pi(A)\right)} + (1-2T) \\
&& + 2T\,\sqrt{\left(\frac 12-\frac{\hat\Phi_T(A)}{2T}\pi(A^c)\right)\left(\frac 12+\frac{\hat\Phi_T(A)}{2T}\pi(A)\right)} \\
&\leq& 2T\sqrt{1-(\hat\Phi_T(A)/2T)^2} + 1-2T \leq 1 - \hat\Phi_T(A)^2/4T
\end{eqnarray*}
The final line was by the relation $\sqrt{XY} + \sqrt{(1-X)(1-Y)} \leq \sqrt{1-(X-Y)^2}$; see Lemma \ref{lem:inequality} in the Appendix for a proof.
If $T=\min\{\alpha,t\}=t$ then the theorem is immediate. If $T=\min\{\alpha,t\}=\alpha$ then use the relation $\hat\Phi_T(A) \geq \frac{\alpha}{t}\,\hat\Phi_t(A)$.
\end{proof}



\section{New Canonical path bounds} \label{sec:alternating-paths}

When ergodic flow leaving a set $A$ is not heavily concentrated at a few vertices then use of threshold conductance methods can improve substantially over use of conductance.
In this section we explore such a situation by considering the case when canonical paths are well distributed among the vertices, and thus the ergodic flow appearing on edges of canonical paths is not heavily concentrated at any vertex.

What is most novel about our approach is that it applies to general finite Markov chains, i.e. non-reversible and non-lazy, a result made possible by using the paths to bound intersection threshold conductance $\overphi_t(A)$.
Ordinary canonical paths cannot be used in this setting as there are walks which do not converge but for which it is easy to construct canonical paths,
such as the simple random walk on a cycle with an even number of vertices (i.e. $\P(i,i-1\mod n)=\P(i,i+1\mod n)=1/2$) and the clockwise walk on a cycle of odd length (i.e. $\P(i,i+1\mod n)=1$).
However, if path length is required to have some {\em parity} (all even or all odd) then the simple walk on an even cycle no longer has natural canonical paths, while a requirement that the path {\em alternate} between edges of $\P$ and $\P^*$ likewise prevents the clockwise walk from having canonical paths.
In fact, these two conditions will be sufficient.

\begin{definition}
A {\em canonical alternating path} $\gamma_{xy}$ is an even length path from $x$ to $y$ alternating between valid transitions of $\P$ and $\P^*$:
\[
x=x_0\xrightarrow{\P}x_1\xrightarrow{\P^*}x_2\xrightarrow{\P}x_3\xrightarrow{\P^*}\cdots\xrightarrow{\P^*} x_{2n}=y
\]
Equivalently, the path alternates between forward and reversed edges of $\P$:
\[
x=x_0\rightarrow x_1 \leftarrow x_2 \rightarrow x_3 \leftarrow \cdots \rightarrow x_{2n-1} \leftarrow x_{2n}=y
\]

Let $\Gamma=\{\gamma_{xy}:\,x,y\in V\}$ be a set including a canonical alternating path for each {\em ordered} pair of distinct vertices $x,y\in V$.

Define $v\in\gamma_{xy}$ if $v=x_{2i+1}$ for some $i$, i.e. $v$ is the terminal point of an edge in $\gamma_{xy}$.

Define $(u,v)\in\gamma_{xy}$ if $u\rightarrow v$ or $v\leftarrow u$ appears in path $\gamma_{xy}$, i.e. $(u,v)$ is an edge of $\P$ or $(v,u)$ an edge of $\P^*$ in the path.
\end{definition}

\noindent
The definitions of $v\in\gamma_{xy}$ and $(u,v)\in\gamma_{xy}$ can be relaxed to allow both initial and terminal endpoints of edges, but this leads to weaker results.

As discussed in the Introduction, different canonical path methods require different types of paths.
However, each of these types of paths induce natural canonical alternating paths, as shown below, so our new definition provides a unifying framework:

\begin{itemize}
\item For a finite reversible walk with holding probability $\alpha>0$ Diaconis and Strook \cite{DS91.1} and Sinclair \cite{Sin92.1} use simple paths:
\[
x=x_0\rightarrow x_1 \rightarrow x_2 \rightarrow x_3 \rightarrow \cdots \rightarrow x_n=y
\]
Such paths can be used to study non-reversible walks as well \cite{MT06.1}.
In either case, to construct an alternating path insert a self-loop at each vertex after the initial one and consider the loop as a transition of $\P^*$ since $\P^*(v,v)=\P(v,v)$:
\[
x=x_0\rightarrow \overloop{x_1} \rightarrow \overloop{x_2} \rightarrow \cdots \rightarrow \overloop{x_n}=y
\]

\item For a general finite reversible walk Diaconis and Strook \cite{DS91.1} use paths in which $\gamma_{xy}$ is constructed for $x\neq y$ with no parity requirement, but with looping paths $\gamma_{xx}$ which must be odd length.
If the length $|\gamma_{xy}|$ is even then it is also an alternating path since $\P^*=\P$ for reversible walks, while if $|\gamma_{xy}|$ is odd then $\gamma_{xy}$ followed by $\gamma_{yy}$ will be an alternating path.

\item For a general finite walk Mihail \cite{Mih89.1} and Fill \cite{Fill91.1} show that it suffices to study conductance of $\P\P^*$ or canonical paths with edges in $\P\P^*$.
A path $\gamma_{xy}$ with edges from $\P\P^*$ can be made into an alternating path by replacing each edge $\P\P^*(x,y)>0$ with a pair of edges $x\xrightarrow{\P}z\xrightarrow{\P^*}y$ such that $\P(x,z)>0$ and $\P^*(z,y)>0$.

\item In an earlier version of this paper we required odd length paths alternating between $\P$ and $\P^*$, including looping paths $\gamma_{xx}$.
Given $x\neq y$ the path $\gamma_{xy}$ followed by $\gamma_{yy}$ in reverse will be an alternating path.
\end{itemize}

Notions of both edge and vertex congestion will be required in our results:

\begin{definition} \label{def:rho_v}
Given (ordinary or alternating) canonical paths $\Gamma$ between every pair of distinct vertices $x,y\in V$ the {\em vertex congestion} is
\[
\rho_v = \rho_v(\Gamma) = \max_{v\in V} \frac{1}{2\pi(v)}\,\sum_{\gamma_{xy}\cup\gamma_{yx}\ni v} \pi(x)\pi(y)
\]
and the {\em edge congestion} is
\[
\rho_e = \rho_e(\Gamma) = \max_{e=(u,v)\in E} \frac{1}{2\pi(u)\P(u,v)}\,\sum_{\gamma_{xy}\cup\gamma_{yx}\ni (u,v)} \pi(x)\pi(y)
\]
\end{definition}

The redefined edge congestion is no larger than the earlier definition and is equivalent to it when, as is typically the case, the path $\gamma_{xy}$ does not use any of the same (directed) edges as $\gamma_{yx}$, e.g. in the reversible case $\gamma_{xy}$ can be assumed to traverse the edges of $\gamma_{yx}$ in reverse.
The conductance and mixing time bounds to be shown in this section are weaker for larger $\rho_e$, so using the earlier definition will give weaker but still valid results.

This brings us to the main result of the section, a lower bound on the intersection threshold conductance $\overphi_t$ in terms of canonical alternating paths.

\begin{lemma} \label{lem:alternating-conductance}
Given canonical alternating paths $\Gamma$ then
\[
\overphi_{\rho_v/\rho_e} \geq 1/2\rho_e
\]
If, instead, $\Gamma$ consists of ordinary canonical paths then
\[
\hat\Phi_{\rho_v/\rho_e} \geq 1/\rho_e
\]
\end{lemma}

The proof is given later.
Combining this with our results on threshold conductances leads to new upper bounds on mixing time:

\begin{theorem} \label{thm:alternating-mixing}
Consider a finite ergodic Markov kernel $\P$ with canonical ordinary or alternating paths.
Assume $\epsilon\leq 1$.
Let
\[
\P_0^*(\Gamma) = \min\left\{\P^*(b,a):\,\exists \gamma_{xy}\ni(a,b)\right\}
\]
be the smallest transition in the reversal $\P^*(b,a)=\frac{\pi(a)\P(a,b)}{\pi(b)}$ of the edges in the paths.

If $\Gamma$ consists of canonical alternating paths then $\P$ has mixing time
\[
\tau(\epsilon)
 \leq 48 \rho_e \max\left\{\rho_v,\,\frac 12\right\} 
                     \left(\log\frac{1}{\epsilon\sqrt{\pi_0}} - \frac 14\left(\log(4\rho_v\rho_e\P_0^*(\Gamma))-1\right)^+\right)
\]
where $x^+=\max\{x,0\}$.

If, instead, $\Gamma$ consists of ordinary canonical paths then
\[
\tau(\epsilon) \leq 4\rho_v\,\max\left\{\frac{\rho_v}{\alpha},\,\rho_e\right\}\,\left(\log\frac{1}{\epsilon\sqrt{\pi_0}} - \frac 14\,\left(\log (\rho_v\rho_e\P_0^*(\Gamma))-1\right)^+\right)
\]
\end{theorem}

\noindent
These can be simplified by using the relations:
\begin{equation} \label{eqn:congestions}
\rho_v \leq \rho_e \leq \frac{\rho_v}{\P_0^*(\Gamma)} < \frac{1}{2\pi_0\P_0^*(\Gamma)}
\end{equation}
With ordinary paths a lower bound of $\rho_v\geq 1-\min\pi(v)$ also holds, while alternating canonical paths have the weaker $\rho_v\geq \frac 12\,(1-\max\pi(v))$.
Although the mixing bounds of the theorem are not monotone, they do hold for any upper bounds on $\rho_v$ and $\rho_e$ which satisfy the constraints just mentioned.

The $\max\left\{\frac{\rho_v}{\alpha},\,\rho_e\right\}$ term is just an upper bound on the edge congestion when alternating paths are constructed by adding self-loops at all but the initial vertices of the ordinary paths. 
This maximum and the correction term with $\log(\rho_v\rho_e\P_0^*(\Gamma))$ are not simply artifacts of the proof.
A simple illustration of this is the walk on a cycle $\ZZ_n$ with $\P(i,i)=\alpha\in(0,1)$ and $\P(i,i+1)=1-\alpha$. With the obvious choice of paths this has $\rho_v=\frac{n-1}{2}$, $\rho_e=\frac{n-1}{2(1-\alpha)}$, $\P_0^*(\Gamma)=1-\alpha$, $\pi_0=1/n$ and mixing time $\tau(\epsilon)=\Theta\left(\frac{n^2}{\min\{\alpha,1-\alpha\}}\log\frac{1}{\epsilon}\right)$.
In contrast, the lazy walk on the complete graph $K_n$ has mixing time $\tau(1/4)=\Omega(\log n)$ and so the $\pi_0$ term is sometimes necessary.

\begin{proof}
First consider ordinary canonical paths.
The theorem is trivial if $|V|=1$, so assume that $|V|\geq 2$ and $0 < \min\pi(v)\leq 1/2$.
Lemma \ref{lem:alternating-conductance} and Theorem \ref{thm:threshold-conductances} bound the root profile as
\[
\psi(A)\geq \frac{\min\{\alpha,\,\rho_v/\rho_e\}}{4(\rho_v/\rho_e)^2}\,\frac{1}{\rho_e^2}
 = \frac{1}{4\rho_v\max\{\rho_v/\alpha,\,\rho_e\}}
\]
Evaluating the integral in the convex case of Theorem \ref{thm:evomixing} gives the result when $\rho_v\rho_e\P_0^*(\Gamma)\leq e$, so assume that $\rho_v\rho_e\P_0^*(\Gamma)> e$. 
Since the walk is ergodic then the state space is strongly connected, and so $\Q(A,A^c)>0$ for every set $A\not\in\{\emptyset,V\}$.
In particular, $\Q_{\P_0^*(\Gamma)}(A,A^c)\geq\pi_0\P_0^*(\Gamma)$, and similarly $\Q_{\P_0^*(\Gamma)}(A^c,A)\geq\pi_0\P_0^*(\Gamma)$, and so $\hat\Phi_{\P_0^*(\Gamma)}(A) \geq \frac{\pi_0\P_0^*(\Gamma)}{\pi(A)\pi(A^c)}$.
By Theorem \ref{thm:threshold-conductances}
\[
\psi(A)\geq \max\left\{
 \frac{1}{4\rho_v\max\{\rho_v/\alpha,\,\rho_e\}},\,
 \frac{\pi_0^2\min\{\alpha,\P_0^*(\Gamma)\}}{4\pi(A)^2}
\right\}
\]
The convexity condition of Theorem \ref{thm:evomixing} is easily verified for this lower bound and so
\[
\tau(\epsilon)\leq \int_0^c \frac{2r\,dr}{\pi_0^2\min\{\alpha,\P_0^*(\Gamma)\}}
  + 2\rho_v\max\left\{\frac{\rho_v}{\alpha},\,\rho_e\right\}\,\int_{c}^{1/\epsilon^2} \frac{dr}{r}
\]
if $c=\pi_0\sqrt{\min\{\alpha,\P_0^*(\Gamma)\}\,\rho_v\max\{\rho_v/\alpha,\rho_e\}} \leq 1/\epsilon^2$.
The theorem follows by integrating and using equations \eqref{eqn:congestions}, $\rho_v\rho_e\P_0^*(\Gamma)\geq 1$ and $1/2 \leq \rho_v < 1/2\pi_0$.

For alternating canonical paths $\psi(A)\geq \left(48\max\{\rho_v,1/2\}\rho_e\right)^{-1}$ by Lemma \ref{lem:alternating-conductance} and Theorem \ref{thm:threshold-conductances}.
Theorem \ref{thm:evomixing} shows the mixing bound when $4\rho_v\rho_e\P_0^*(\Gamma)\leq e$.
To improve $\psi(A)$ use
$
\overphi_{\P_0^*(\Gamma)} \geq \frac{\pi_0\P_0^*(\Gamma)}{\pi(A)}
$
because there is at least one alternating path between some $x\in A$ and $y\in A^c$, and this path will contain some vertex $v\in\gamma_{xy}$ with incoming edges from both $A$ and $A^c$.
Again use Theorem \ref{thm:threshold-conductances} to lower bound $\psi(A)$, and this time split the mixing time integral of Theorem \ref{thm:evomixing} at $c=2\pi_0\sqrt{\rho_v\rho_e\P_0^*(\Gamma)}$.
\end{proof}

We return now to the proof of our main result:

\begin{proof}[Proof of Lemma \ref{lem:alternating-conductance}]
Consider either the ordinary canonical path or alternating path case.
The stationary distribution and edge capacity can be lower bounded using paths:
\begin{eqnarray*}
\forall v\in V:& \displaystyle
\pi(v) &\geq \frac{1}{2\rho_v}\,\sum_{\gamma_{xy}\cup \gamma_{yx}\ni v} \pi(x)\pi(y) \\
    &&\geq \frac{1}{2\rho_v}\,\sum_{\substack{(x,y):\,\exists u\in A,\,\\(u,v)\in\gamma_{xy}\cup\gamma_{yx}}} \pi(x)\pi(y) \\
\forall u\in A,\,v\in V:& \displaystyle
\pi(u)\P(u,v) &\geq \frac{1}{2\rho_e}\,\sum_{\gamma_{xy}\cup\gamma_{yx}\ni(u,v)} \pi(x)\pi(y)\,.
\end{eqnarray*}
Hence, for any $v\in V$:
\begin{eqnarray}
\Q_{\rho_v/\rho_e}(A,v)
 &=&  \min\left\{\sum_{u\in A} \pi(u)\P(u,v),\,\frac{\rho_v}{\rho_e}\,\pi(v)\right\} \nonumber \\
 &\geq& \frac{1}{2\rho_e}\,\min\left\{ \sum_{u\in A} \sum_{\gamma_{xy}\cup\gamma_{yx}\ni(u,v)} \pi(x)\pi(y),\,
                               \sum_{\substack{(x,y):\,\exists u\in A,\,\\(u,v)\in\gamma_{xy}\cup\gamma_{yx}}} \pi(x)\pi(y)
                               \right\} \nonumber \\
 &=&  \frac{1}{2\rho_e}\,\sum_{\substack{(x,y):\,\exists u\in A,\,\\(u,v)\in\gamma_{xy}\cup\gamma_{yx}}} \pi(x)\pi(y) \label{eqn:t-flow}
\end{eqnarray}

In the ordinary canonical path case it follows from \eqref{eqn:t-flow} that
\begin{eqnarray*}
\lefteqn{ \Q_{\rho_v/\rho_e}(A,A^c) = \sum_{v\in A^c} \Q_{\rho_v/\rho_e}(A,v) } \\
 &\geq& \sum_{v\in A^c} \frac{1}{2\rho_e}\,\sum_{\substack{(x,y):\,\exists u\in A,\,\\(u,v)\in\gamma_{xy}\cup\gamma_{yx}}} \pi(x)\pi(y) \\
 &\geq& \frac{1}{2\rho_e}\,\sum_{(x,y)\in A\times A^c} 2\pi(x)\pi(y) 
  =     \frac{\pi(A)\pi(A^c)}{\rho_e}
\end{eqnarray*}
The second inequality is because a path from some $x_0\in A$ to $x_n\in A^c$ must have some $x_i\in A$ and $x_{i+1}\in A^c$.
Taking $A \leftarrow A^c$ shows that $\Q_{\rho_v/\rho_e}(A^c,A) \geq \frac{\pi(A)\pi(A^c)}{\rho_e}$ as well.
It follows that $\hat\Phi_{\rho_v/\rho_e}(A)\geq 1/\rho_e$.

If the paths are alternating then \eqref{eqn:t-flow} shows that
\begin{eqnarray*}
\lefteqn{ \sum_{v\in V} \min\{\Q_{\rho_v/\rho_e}(A,v),\,\Q_{\rho_v/\rho_e}(A^c,v)\} } \\
 &\geq& \sum_{v\in V} \frac{1}{2\rho_e}\,\sum_{\substack{(x,y): \exists u\in A,\,w\in A^c:\\ (u,v)\in\gamma_{xy}\cup\gamma_{yx},\\ (w,v)\in\gamma_{xy}\cup\gamma_{yx}}} \pi(x)\pi(y) \\
 &\geq& \frac{1}{2\rho_e}\,\sum_{(x,y)\in A\times A^c} 2\pi(x)\pi(y) 
  \geq \frac{\pi(A)\pi(A^c)}{\rho_e}
\end{eqnarray*}
The second inequality is because an (even length) alternating path from $x_0\in A$ to $x_{2n}\in A^c$ must have some $x_{2i}\in A$ and $x_{2(i+1)}\in A^c$, and so is counted when $v=x_{2i+1}$.
It follows that $\overphi_{\rho_v/\rho_e}(A)\geq 1/2\rho_e$.
\end{proof}


\begin{remark} \label{rmk:averages}
Ordinary canonical paths can be applied to three types of mixing bounds for lazy walks, with lead term $\rho_e^2$ for Jerrum and Sinclair's result, $\rho_e\ell$ for the Poincar\'e bound, and $\rho_e\rho_v$ for ours.
Since $\rho_v\leq\rho_e$ then our bound improves on Jerrum and Sinclair's.
To compare to the Poincar\'e bound define average vertex congestion and average path length by
\[
\bar{\rho}_v = \sum_{v\in V} \pi(v) \left[ \frac{1}{2\pi(v)} \sum_{\gamma_{xy}\cup\gamma_{yx}\ni v} \pi(x)\,\pi(y) \right]
\qquad \textrm{and}\qquad
\bar{\ell} = \frac{\sum_{x\neq y} \pi(x)\,\pi(y)\,|\gamma_{xy}|}{\sum_{x\neq y} \pi(x)\pi(y)} 
\] 
respectively. Then
\begin{eqnarray*}
\bar{\rho}_v =& \displaystyle \frac 12\,\sum_{x\neq y} \sum_{v\in \gamma_{xy}\cup\gamma_{yx}} \pi(x)\,\pi(y) 
    &\leq \sum_{x\neq y} \pi(x)\,\pi(y)\,|\gamma_{xy}| \nonumber \\
    =& \displaystyle \bar{\ell}\,\sum_{x\neq y} \pi(x)\,\pi(y)
    &= \bar{\ell}(1-\|\pi\|_2^2)\ .
\end{eqnarray*}
Likewise $\bar{\rho}_v \geq \frac{\bar{\ell}+1}{2}(1-\|\pi\|_2^2)$, so $\bar{\rho}_v = \Theta\left(\bar{\ell}(1-\|\pi\|_2^2\right)$.
When a few paths are very long or the distribution is concentrated near a single vertex then it is likely that $\rho_v\ll\ell$ and our bound is best.
However, it is more often the case that path length varies little and a bottleneck causes a few states to have high vertex congestion, and so $\ell\ll\rho_v$.
In contrast to the lazy case, when the holding probability is small then our result can be significantly better, even when $\ell\ll\rho_v$; see the next section for examples.
\end{remark}


\section{Examples} \label{sec:examples}

To demonstrate our method we give two examples where the new canonical path theorems extend previously known bounds into the general non-reversible non-lazy case. First, the classical problem of the max-degree walk on a graph, then the more interesting case of walks on Cayley graphs, i.e. random walks on groups.
 
\begin{example}
An Eulerian multigraph is a strongly connected graph with in-degree=out-degree at each vertex, a natural generalization of the undirected multigraph into the directed graph setting.
Suppose an Eulerian multigraph has $n$ vertices and maximum out-degree $d$.
Let $d(x,y)$ denote the number of directed edges from $x$ to $y$, so that $d(x)=\sum_y d(x,y)$ is the out-degree of $x$.
The max-degree walk has $\P(x,y)=\frac{d(x,y)}{d}$ if $y\neq x$ and $\P(x,x)=1-\frac{d(x)-d(x,x)}{d}$.
The stationary distribution $\pi=1/n$ is uniform and the walk is lazy if $d(x)-d(x,x)\leq d/2$ at each vertex, i.e. there are many self-loops.

Suppose that every vertex has a self-loop. Then holding probability $\alpha\geq 1/d$, which is sufficient to prevent significant periodicity effects.
For every $x\neq y$ let $\gamma_{xy}$ be a path from $x$ to $y$.
The congestions are at worst
\[
\rho_v\leq\frac{\sum_{x\neq y}\pi(x)\pi(y)}{2\pi_0}\leq\frac{n}{2}
\quad\textrm{and}\quad
\rho_e\leq\frac{\sum_{x\neq y}\pi(x)\pi(y)}{2\pi_0\P_0(\Gamma)} \leq \frac{dn}{2}
\]
and so by Theorem \ref{thm:alternating-mixing}
\[
\tau(\epsilon)\leq dn^2\,\log\frac{2}{\epsilon}
\]

Without the self-loop requirement, if the graph is connected under canonical alternating paths then it is still true that $\rho_v\leq\frac n2$ and $\rho_e\leq \frac{nd}{2}$ and so by Theorem \ref{thm:alternating-mixing}
\[
\tau(\epsilon)\leq 12\,dn^2\log\frac{2}{\epsilon}\,.
\]
Recall that strong connectivity is not enough for ergodicity, e.g. the cycle walk $\P(i,i+1\mod n)=1$ on $\ZZ_n$.

In contrast, the Poincar\'e method introduces an extra factor of $\alpha^{-1}$, and so it can match the self-looping case only for a walk with constant (in $n$ and $d$) holding probability.
Diaconis and Strook's extension to non-lazy walks works only for the reversible case (i.e. $d(x,y)=d(y,x)$ for every $x,y$).
Mihail and Fill's extension with $\P\P^*$ replaces the $(\min_{\P(x,y)>0} \P(x,y))^{-1}$ term in $\rho_e$ with $(\min_{\P\P^*(x,y)>0} \P\P^*(x,y))^{-1}$, which typically replaces the order $d$ term with order $d^2$.
\end{example}

\begin{example}\label{ex:Cayley}
The Cayley graph of a group $G$ with (non-symmetric) generating set $S\subset G-\{id\}$, i.e. $\bigcup_{n=0}^\infty S^n=G$, 
has edge set $(g,gs)$ for all $g\in G$, $s\in S$.
If $p:\,G\to[0,1]$ is a probability distribution supported on $S\cup\{id\}$ then $\P(g,gs)=p(s)$ defines a Markov chain with uniform stationary distribution $\pi=1/|G|$.
Represent each $g\in G$ as a product of generators $g=s_1\,s_2\,\cdots\,s_k$,
define $\Delta=\max |g|$ to be the length of the longest such representation,
and let $N(g,s)\leq\Delta$ denote the number of times generator $s$ appears in the representation of $g$.

Babai \cite{Bab91.1} showed $\tau(\epsilon)=O(\frac{\Delta^2}{\min_{s\in S}p(s)}\log \frac{|G|}{\epsilon})$ for the lazy walk with symmetric generating set, i.e. $p(id)\geq 1/2$ with $S=S^{-1}$ and $\forall s\in S:\,p(s)=p(s^{-1})$.
Diaconis and Saloff-Coste \cite{DSC93.2} use (ordinary) canonical paths to bound the spectral gap by the Poincar\'e approach. This can be plugged into spectral gap bounds on mixing time (e.g. Corollary 2.15 of \cite{MT06.1}), leading to the following generalizations of Babai's result to the symmetric and non-symmetric cases respectively:
\begin{eqnarray}
\tau(\epsilon) &\leq& \max\left\{\frac{1}{2\,p(id)},\,\Delta\max_{g\in G,\,s\in S} \frac{N(g,s)}{p(s)}\right\}\,\left(\frac 12\log |G| + \log\frac{1}{\epsilon}\right) \nonumber \\
\tau(\epsilon) &\leq& \frac{\Delta}{p(id)}\,\max_{g\in G,\,s\in S} \frac{N(g,s)}{p(s)}\,\left(\frac 12\log |G| + \log\frac{1}{\epsilon}\right) \label{eqn:cayley-old}
\end{eqnarray}

Consider now our new method of canonical alternating paths.
Let $\Delta_{alt}$ be the diameter measured using canonical alternating paths, i.e. $\Delta_{alt}=2\,\min\{N:\,G=\bigcup_{n=0}^N (SS^{-1})^n\}$, and let $N_{alt}(g,s)$ count frequency of generators when group elements are written in terms of alternating canonical paths, e.g. $g=s_1\,s_2^{-1}\,s_3\,s_4^{-1}\,\cdots\,s_{2n-1}\,s_{2n}^{-1}$.
Congestions for canonical alternating paths lead to the bounds
\[
\rho_v < \frac{\Delta_{alt}}{4},\qquad
\rho_e < \frac 12\,\max_{g\in G,\,s\in S}\,\frac{N_{alt}(g,s)}{p(s)}
\]
The proof is left to the Appendix as it uses essentially the same approach used for ordinary canonical paths in \cite{DSC93.2}.
Theorem \ref{thm:alternating-mixing} leads to a mixing bound of
\begin{equation} \label{eqn:Cayley-alternating}
\tau(\epsilon)
 \leq 6\,\Delta_{alt}\,\max_{g\in G,\,s\in S}\,\frac{N_{alt}(g,s)}{p(s)}\left(\frac 12\,\log\frac{4|G|}{\Delta_{alt}}+\log\frac{1}{\epsilon}\right)
\end{equation}
To see that this generalizes Diaconis and Saloff-Coste's non-symmetric result recall the construction of alternating paths from ordinary canonical paths by adding a self-loop at each vertex except the starting point.
Such a set of alternating paths will have $\Delta_{alt}=2\Delta$, while the self-loops cause
\[
\max_{g\in G,\,s\in S}\,\frac{N_{alt}(g,s)}{p(s)} = \max\left\{\frac{\Delta}{p(id)},\,\max_{g\in G,\,s\in S}\,\frac{N(g,s)}{p(s)}\right\}
\leq \frac{1}{p(id)}\,\max_{g\in G,\,s\in S}\,\frac{N(g,s)}{p(s)}
\]
The inequality was because 
$\Delta = \max_{g\in G} \sum_{s\in S} p(s)\,\frac{N(g,s)}{p(s)} \leq \max_{g\in G}\max_{s\in S}\frac{N(g,s)}{p(s)}$.

It follows that, up to a constant, if ordinary canonical paths are replaced by canonical alternating paths then past results for walks on Cayley graphs hold even when there is no holding probability.
\end{example}


\bibliographystyle{plain}
\bibliography{../references}

\begin{thebibliography}{10}

\bibitem{Bab91.1}
L.~Babai.
\newblock Local expansion of vertex-transitive graphs and random generation in
  finite groups.
\newblock {\em Proceedings of the 23rd Annual ACM Symposium on Theory of
  Computing (STOC 1991)}, pages 164--174, 1991.

\bibitem{DF90.1}
P.~Diaconis and J.~Fill.
\newblock Strong stationary times via a new form of duality.
\newblock {\em The Annals of Probability}, 18(4):1483--1522, 1990.

\bibitem{DSC93.2}
P.~Diaconis and L.~Saloff-Coste.
\newblock Comparison techniques for random walk on finite groups.
\newblock {\em The Annals of Probability}, 21(4):2131--2156, 1993.

\bibitem{DS91.1}
P.~Diaconis and D.~Stroock.
\newblock Geometric bounds for eigenvalues of markov chains.
\newblock {\em The Annals of Applied Probability}, 1:36--61, 1991.

\bibitem{Fill91.1}
J.~Fill.
\newblock Eigenvalue bounds on convergence to stationarity for nonreversible
  markov chains, with an application to the exclusion process.
\newblock {\em The Annals of Applied Probability}, 1(1):62--87, 1991.

\bibitem{JS88.1}
M.~Jerrum and A.~Sinclair.
\newblock Conductance and the rapid mixing property for markov chains: the
  approximation of the permanent resolved.
\newblock {\em Proceedings of the 20th Annual ACM Symposium on Theory of
  Computing (STOC 1988)}, pages 235--243, 1988.

\bibitem{KLM06.1}
R.~Kannan, L.~Lov\'asz, and R.~Montenegro.
\newblock Blocking conductance and mixing in random walks.
\newblock {\em Combinatorics, Probability and Computing}, 15(4):541--570, 2006.

\bibitem{LS88.1}
G.~Lawler and A.~Sokal.
\newblock Bounds on the $l^2$ spectrum for markov chains and markov processes:
  a generalization of cheeger's inequality.
\newblock {\em Transactions of the American Mathematical Society},
  309:557--580, 1988.

\bibitem{Mih89.1}
M.~Mihail.
\newblock Conductance and convergence of markov chains-a combinatorial
  treatment of expanders.
\newblock {\em 30th Annual Symposium on Foundations of Computer Science}, pages
  526--531, 1989.

\bibitem{MT06.1}
R.~Montenegro and P.~Tetali.
\newblock {\em Mathematical Aspects of Mixing Times in Markov Chains}, volume
  1:3 of {\em Foundations and Trends in Theoretical Computer Science}.
\newblock NOW Publishers, Boston-Delft, June 2006.

\bibitem{MP05.1}
B.~Morris and Y.~Peres.
\newblock Evolving sets, mixing and heat kernel bounds.
\newblock {\em Probability Theory and Related Fields}, 133(2):245--266, 2005.

\bibitem{Sin92.1}
A.~Sinclair.
\newblock Improved bounds for mixing rates of markov chains and multicommodity
  flow.
\newblock {\em Combinatorics, Probability and Computing}, 1(4):351--370, 1992.

\end{thebibliography}


\appendix
\section{Appendix}

The following inequality was used in the proof of Theorem \ref{thm:threshold-conductances}.

\begin{lemma} \label{lem:inequality}
If $X,\,Y\in[0,1]$ then
\[
g(X,Y)=\sqrt{X\,Y}+\sqrt{(1-X)(1-Y)} \leq \sqrt{1-(X-Y)^2}\,.
\]
\end{lemma}

\begin{proof}
Observe that
\[
g(X,Y)^2 = 1-(X+Y)+2\,X\,Y + \sqrt{[1-(X+Y)+2\,X\,Y]^2 - [1-2(X+Y)+(X+Y)^2]}\,.
\]

Now, $\sqrt{A^2-B} \leq A-B$ if $A^2\geq B$, $A\leq\frac{1+B}{2}$ and $A\geq B$ (square both sides to show this). 
These conditions are easily verified with $A=1-(X+Y)+2\,X\,Y$ and $B=1-2(X+Y)+(X+Y)^2$, and so
\begin{eqnarray*}
g(X,Y)^2 &\leq& 2\left[1-(X+Y)+2\,X\,Y\right] - \left[1-2(X+Y)+(X+Y)^2\right] \\
    &=&  1+2\,X\,Y-X^2-Y^2 = 1-(X-Y)^2
\end{eqnarray*}
\end{proof}

In order to study walks on Cayley graphs it was claimed that congestion bounds of Diaconis and Saloff-Coste \cite{DSC93.2} generalize easily. We show this here.

\begin{lemma} \label{lem:groups}
Consider group $G$ with (non-symmetric) generating set $S$.
Use the notation of Example \ref{ex:Cayley} to describe a walk on the Cayley graph of $G=\langle S\rangle$.

There are ordinary canonical paths with
\[
\rho_v < \Delta,\qquad
\rho_e < \max_{g\in G,\,s\in S} \frac{N(g,s)}{p(s)}\,.
\]
If $p$ is symmetric, i.e. $\forall s\in S:\,p(s)=p(s^{-1})$, then $\rho_v < \frac{\Delta+1}{2}$.

There are canonical alternating paths with
\[
\rho_v < \frac{\Delta_{alt}}{4},\qquad
\rho_e < \frac 12\,\max_{g\in G,\,s\in S}\,\frac{N_{alt}(g,s)}{p(s)}
\]
and $N_{alt}(g,s)\leq \frac 12\,\Delta_{alt} < |G|$ for every $g\in G,\,s\in S$.
\end{lemma}

\begin{proof}
First consider ordinary canonical paths.
Given $x,y\in G$ let $g=x^{-1}y=s_1\,s_2\,\cdots\,s_k$ and define path $\gamma_{x,y}$ by $x\to xs_1\to\cdots\to xg=y$.
Recall that $\pi=1/|G|$ is uniform.

To bound vertex-congestion observe that the same number of paths pass through each vertex, because if $\gamma_{x,y}$ includes vertex $v$ then $\gamma_{wv^{-1}x,\,wv^{-1}y}$ includes vertex $w$, and vice-versa. Then
\begin{eqnarray}
\rho_v = \bar{\rho}_v
 &\leq& \frac{1}{|G|}\,\sum_{g\in G} \frac{1}{2\pi(g)}\,\sum_{(x,y):\,g\in \gamma_{xy}\cup\gamma_{yx}} \pi(x)\pi(y) \nonumber \\
 &\leq& \frac 12\,\sum_{x\neq y} \pi(x)\pi(y)\,(|\gamma_{xy}|+|\gamma_{yx}|) \label{eqn:paths-ordinary} \\
 &\leq& \Delta\,\left(1-\sum_{g\in G}\pi(g)^2\right) = \Delta\left(1 - \frac{1}{|G|}\right) \nonumber
\end{eqnarray}

When $p$ is symmetric then assume the representation for $g^{-1}$ to be the inverse of that for $g$, i.e. if $g=s_1\,s_2\,\cdots\,s_k$ then $g^{-1}=s_k^{-1}\,s_{k-1}^{-1}\,\cdots\,s_1^{-1}$.
This does not increase $\Delta$ so it can only improve the bound on $\rho_v$.
If $g\in\gamma_{xy}$ then $g\in\gamma_{yx}$, and vice-versa, so \eqref{eqn:paths-ordinary} improves to $\frac 12\,\sum_{x\neq y} \pi(x)\pi(y)\,(|\gamma_{xy}|+1)$ and $\rho_v \leq \frac{\Delta+1}{2}\left(1 - \frac{1}{|G|}\right)$.

Now consider edge-congestion.
Without loss assume that $id$ does not appear in any paths.
If $\gamma_{x,y}$ includes edge $(v,vs)$ then $\gamma_{wv^{-1}x,wv^{-1}y}$ includes edge $(w,ws)$, and vice-versa,
and so for fixed $s\in S$ the number of paths through edge $(g,gs)$ is independent of the choice of $g\in G$.
Hence, 
\begin{eqnarray*}
\rho_e &\leq& \max_{s\in S} \frac{1}{|G|}\sum_{g\in G} \frac{1}{2\pi(g)\P(g,gs)}\,
                              \sum_{(x,y):\,(g,gs)\in\gamma_{xy}\cup\gamma_{yx}} \pi(x)\pi(y) \\
 &\leq& \max_{s\in S} \frac{1}{2p(s)}\,2\sum_{x\neq y} N(x^{-1}y,s)\pi(x)\pi(y) \\
 &\leq& \max_{s\in S} \frac{1}{p(s)}\,\max_{g\in G} N(g,s)\left(1 - \frac{1}{|G|}\right)
\end{eqnarray*}

Finally, when alternating canonical paths are used then again assume the representation of each $g^{-1}$ to be the inverse of that for $g$, so that $g\in\gamma_{xy}\Leftrightarrow g\in\gamma_{yx}$ and $(g,h)\in\gamma_{xy}\Leftrightarrow (g,h)\in\gamma_{yx}$. Then, arguing as before,
\begin{eqnarray*}
\rho_v = \bar{\rho}_v
 &=& \frac{1}{|G|}\,\sum_{g\in G} \frac{1}{2\pi(g)}\,\sum_{(x,y):\,g\in\gamma_{xy}\cup\gamma_{yx}} \pi(x)\pi(y) \\
 &=& \frac 12\,\sum_{x\neq y} \pi(x)\pi(y)\,\frac{|\gamma_{xy}|}{2} \\
 &\leq& \frac{\Delta_{alt}}{4}\,\left(1-\sum_{g\in G}\pi(g)^2\right) = \frac{\Delta_{alt}}{4}\,\left(1 - \frac{1}{|G|}\right)
\end{eqnarray*}
Similarly minor changes show that $\rho_e < \frac 12\,\max_{s\in S}\,\frac{1}{p(s)} \,\max_{g\in G} N_{alt}(g,s)$.

For the final statement, in the representation $g=s_1\,s_2^{-1}\,\cdots\,s_{2k-1}\,s_{2k}^{-1}$ remove all even length subcycles, reducing the problem to the case where there are no even length subcycles.
In particular, if $i,j\leq k$ then $s_1\,s_2^{-1}\,\cdots\,s_{2i-1}\,s_{2i}^{-1} = s_1\,s_2^{-1}\,\cdots\,s_{2j-1}\,s_{2j}^{-1} \Rightarrow i=j$.
The guarantees that $\left\{id,\,s_1\,s_2^{-1},\,\ldots\,,\,s_1\,s_2^{-1}\,\cdots\,s_{2k-1}\,s_{2k}\right\}$ is a set of $k+1$ distinct elements, and so $k+1\leq |G|$,
while it also guarantees that $s\,s^{-1}$ never appears and so also $N_{alt}(g,s)\leq|g|/2\leq|G|-1$.
\end{proof}

\end{document}